\documentclass[a4paper,11pt]{article}
\usepackage[multidot]{grffile}
\usepackage{latexsym,amssymb,enumerate,amsmath,epsfig,amsthm}
\usepackage[margin=1in]{geometry}
\usepackage{setspace,color}
\usepackage{graphics}
\usepackage{graphicx}
\usepackage[ruled]{algorithm2e}
\usepackage{empheq}
\usepackage{bm,amsmath}
\usepackage{subcaption}
\usepackage{xcolor}

\usepackage{multirow}
\newcommand\fixup{\kern-\fontcharic\scriptfont2`\"}

\newcommand{\Yuki}[1]{\textcolor{violet}{#1}}

\newcommand{\correct}[1]{{#1}}

\title{Spherical Essentially Non-Oscillatory (SENO) Interpolation}
\author{
Ki Wai Fong\thanks{Department of Mathematics, the Hong Kong University of Science and Technology, Clear Water Bay, Hong Kong. Email: {\bf kwfongab@connect.ust.hk}}
\and
Shingyu Leung\thanks{Department of Mathematics, the Hong Kong University of Science and Technology, Clear Water Bay, Hong Kong. Email: {\bf masyleung@ust.hk}}
}

\date{}

\begin{document}
\thispagestyle{plain}
\maketitle

\begin{abstract}
We develop two new ideas for interpolation on $\mathbb{S}^2$. In this first part, we will introduce a simple interpolation method named \textit{Spherical Interpolation of orDER} $n$ (SIDER-$n$) that gives a $C^{n}$ interpolant given $n \geq 2$. The idea generalizes the construction of the B\'{e}zier curves developed for $\mathbb{R}$. The second part incorporates the ENO philosophy and develops a new \textit{Spherical Essentially Non-Oscillatory} (SENO) interpolation method. When the underlying curve on $\mathbb{S}^2$ has kinks or sharp discontinuity in the higher derivatives, our proposed approach can reduce spurious oscillations in the high-order reconstruction. We will give multiple examples to demonstrate the accuracy and effectiveness of the proposed approaches.
\end{abstract}

\section{Introduction}

We consider the following interpolation problem on $\mathbb{S}^2$. \correct{Given a set of ordered data points $\{\mathbf{p}_i: \|\mathbf{p}_i\|_2=1 \mbox{ for } i=0,\cdots,n\}$ on the unit sphere $\mathbb{S}^2$, we aim to determine a parametrized curve $\{\mathbf{p}(t) \in \mathbb{S}^2 : 0\le t\le n\}$ such that $\mathbf{p}(t=i)=\mathbf{p}_i$ at a set of uniformly spaced $t=i$ for $i=0,\cdots,n$.} This interpolation problem has many other important science and engineering applications \cite{kui02}. One of the most popular applications is on computer graphics when we use a point on $\mathbb{S}^2$ to represent rotations of a rigid body about an arbitrary axis. The interpolation problem can be regarded as the in-between orientation of orientations. Another possible interpretation is an interpolation problem in the special orthogonal group $SO(3)$ \cite{shi09}, which applies to the path planning of a rigid body. We could also find applications of this problem in the quantum field theory from quantum mechanics \cite{adl86}, modeling protein structures \cite{pro14}, molecular dynamics simulation \cite{rap85}, theory in fluid mechanics \cite{ghkr06}, fluid flow visualizations \cite{hanma95}, computations of flexible filaments and fibres in complex fluids \cite{tscfro20,stwk21}, and differential equations \cite{kouxia18} and some applications to dynamics of rigid-bodies \cite{weiterfed06,wil09,udwsch10}.

Several methods exist to interpolate the action as the data points on the unit sphere $\mathbb{S}^2$. For example, based on the quaternion representation \cite{shoemake_85,watwat92,zha97,muk02,barhasben04} of data points on a unit sphere, one has the spherical linear interpolation (SLERP) and the spherical quadrangle interpolation (SQUAD). These two are the most popular and commonly used interpolation methods on the unit sphere. The SLERP interpolation provides the piecewise linear interpolation in geodesic on the unit sphere, while the SQUAD gives a smooth and slightly higher-order reconstruction of the data points. To our understanding, there is no systematic construction of interpolation schemes that provides any order higher than SQUAD. 

There are multiple possible reasons for the absence of high order interpolation scheme. One possible explanation is that data might contain noise in most real-life applications, such as computer graphics or unmanned aerial vehicle (UAV) trajectory planning. There is no significant reason to treat all data very precisely. Therefore, most applications rely on B\'{e}zier curves for smooth curve constructions. In some other applications, however, when we need to solve some related differential equation where the numerical precision of the solution is essential, we might worry about the continuity of the underlying unknown curve and avoid high-order reconstruction. This consideration is common in the usual Euclidean space. If the underlying function is not smooth enough, we could experience the so-called Gibbs phenomenon and obtain spurious oscillations from high-order polynomial reconstructions. Instead of treating all data points equally in the small neighborhood, the essentially non-oscillatory (ENO) scheme \cite{harengoshcha87,shuosh88,shuosh89,shu90} obtains a biased choice of grid stencil that yields the least variations in the polynomial reconstruction. This interpolation method becomes the essential tool in the numerical methods for solving Hamilton-Jacobi equations \cite{jiapen00,zhashu03,serqia06} in the level set applications \cite{oshset88,set99,oshfed03} where the solution might develop kinks or for approximating the hyperbolic conservation laws \cite{shuosh88,shuosh89} when the solution profile might form shocks and discontinuities.

This paper has two primary purposes. We first propose a systematic approach to construct some high-order interpolants from given data points on $\mathbb{S}^2$. We name our approach the \textit{Spherical Interpolation of orDER} $n$ (SIDER-$n$) that produces a $C^{n}$ interpolant given $n \geq 2$. The idea follows the construction of the B\'{e}zier curves based on the composition of multiple SLERPs. Therefore, the formulation is familiar with what has been widely used in the community. However, unlike the typical B\'{e}zier curves, our proposed {approach determines} the control points {in B\'{e}zier} curves to enforce that the reconstruction passes through all given data points. Once we have these high-order reconstructions on $\mathbb{S}^2$, we follow the ENO philosophy and propose a new \textit{Spherical Essentially Non-Oscillatory} (SENO) interpolation method. This approach can provide a smooth, high-order interpolant even when the underlying curve has kinks. This property will be necessary when we aim to develop high-order numerical schemes for solving any quaternion-related differential equations.

The paper is organized as follows. In Section \ref{Sec:Quaternions}, we briefly introduce the quaternion notation and provide the typical interpolation approaches, including SLERP and SQUAD. To increase the regularity of the interpolant, we develop the SIDER and the SENO interpolation method in Section \ref{Sec:Proposal}. Section \ref{Sec:Example} shows various numerical examples to demonstrate the accuracy and effectiveness of the proposed interpolation method.

\section{Quaternions and Interpolations on the Unit Sphere}
\label{Sec:Quaternions}

This section will first introduce the quaternion notation and summarize some essential properties. Based on this quaternion representation, we will briefly introduce the SLERP and SQUAD interpolations. 

Hamilton introduced quaternions to describe rotations and scalings in mid nineteenth century \cite{ham63}. These numbers consist of four dimensions, one real part and a three-dimensional analogy to {the imaginary part of} complex numbers. A quaternion can be written in many forms: \[\underset{\text{real}}{\boxed{a}} + \underset{\text{imaginary}}{\boxed{b \mathbf{i} + c\mathbf{j} + d\mathbf{k}}} = (a, b, c, d) = (\underset{\text{scalar}}{\boxed{a}}, \underset{\text{vector}}{\boxed{\mathbf{u}}}),\] where \(a, b, c, d \in \mathbb{R}\), $\mathbf{u} = (b, c, d) \in \mathbb{R}^3$. The notations $\mathbf{i}$, $\mathbf{j}$ and $\mathbf{k}$ are extensions of {the imaginary part of} complex numbers with the properties that $\mathbf{i}^2=\mathbf{j}^2=\mathbf{k}^2{=\mathbf{ijk}}=-1$, $\mathbf{ij}=\mathbf{k}$ with the {bicyclic permutation with respect to $\mathbf{i}$ that $1\rightarrow \mathbf{i}\rightarrow -1\rightarrow -\mathbf{i}$ and $\mathbf{j}\rightarrow -\mathbf{k}\rightarrow -\mathbf{j}\rightarrow \mathbf{k}$}. Some important properties about quaternions include \begin{itemize}
\item Hamilton product
$$
(a_1, \mathbf{u_1})(a_2, \mathbf{u_2}) = (a_1a_2 - \mathbf{u_1} \cdot \mathbf{u_2}, a_1\mathbf{u_2} + a_2\mathbf{u_1} + \mathbf{u_1} \times \mathbf{u_2}) \, ,
$$
where the notation {$\cdot$ and $\times$ denotes the typical dot and cross product.}
\item \correct{Inverse map
$$\mathbf{q}^{-1} = (a, -\mathbf{u})/(a^2 + b^2 + c^2 + d^2)$$ If $\mathbf{q} = (a, \mathbf{u})$. In particular, if $\mathbf{q}$ is a unit quaternion, $\mathbf{q}^{-1} = (a, -\mathbf{u})$.}
\item Exponential map 
$$
\exp(a, \mathbf{u}) = \exp(a)(\cos \lVert\mathbf{u}\rVert, ((\sin \lVert\mathbf{u}\rVert)/\lVert\mathbf{u}\rVert) \mathbf{u}) \, .
$$ {where the norm notation $\lVert \cdot \rVert$ denotes the 2-norm in this paper, unless otherwise specified.}
\item Logarithm map
$$
\ln(a, \mathbf{u}) = \left(\ln \sqrt{a^2 + \lVert \mathbf{u} \rVert^2}, \frac{1}{\lVert \mathbf{u} \rVert} 
\arccos \left( \frac{a}{\sqrt{a^2 + \lVert \mathbf{u} \rVert^2}} \right) 
\mathbf{u} \right) \, .
$$ 
\item Power map 
\begin{align*}
&\hspace{1.33em} (a, \mathbf{u})^{f(t)} = \exp(f(t)\ln(a, \mathbf{u})) \\
&= \exp(f(t)\ln \sqrt{a^2 + \lVert \mathbf{u} \rVert^2}, (f(t)k/\lVert \mathbf{u} \rVert) \mathbf{u}) \\
&= \left((a^2 + \lVert \mathbf{u} \rVert^2)^{f(t)/2}\cos (f(t)k), (a^2 + \lVert \mathbf{u} \rVert^2)^{f(t)/2} \left[\sin (f(t)k)/\lVert\mathbf{u}\rVert\right] \mathbf{u}\right) \, ,
\end{align*}
where $k = \arccos \left(a/\sqrt{a^2 + \lVert \mathbf{u} \rVert^2}\right)$. When a quaternion has its {2-}norm $\sqrt{a^2 + b^2 + c^2 + d^2}$ equal to one, we call them a unit quaternion. If $(a, \mathbf{u})$ is a unit quaternion, then 
$$
(a, \mathbf{u})^{f(t)} = \left(\cos (f(t)k), \left[ \frac{\sin (f(t)k)}{\lVert \mathbf{u} \rVert} \right] \mathbf{u} \right) \, .
$$
\end{itemize}

Because we can use these unit quaternions to define rotation, we also call these quaternions rotation quaternions. With proper definitions, they can rotate a position vector defined in either $\mathbb{S}^2$ or $\mathbb{R}^3$ while preserving the length of the vector. To see this, we express a unit quaternion as 
$$
(a, b, c, d) = (a, \mathbf{u}) = (\cos (\theta/2), \correct{\sin (\theta/2)}\mathbf{v}) 
$$ 
where $\mathbf{v}$ is a unit vector representing the 3D rotation axis, and $\theta$ is the anticlockwise\correct{/counterclockwise} rotation angle around $\mathbf{v}$ carried by the rotation quaternion. If we want to rotate $\mathbf{p_a}\in \mathbb{S}^2$ with a rotation quaternion $\mathbf{r_{ab}} = (\cos (\theta_{ab}/2), \sin (\theta_{ab}/2)\mathbf{a_{ab}})$, we can first convert $\mathbf{\mathbf{p_a}}$ to a unit quaternion given by $\mathbf{\mathbf{q_a}} = (0, \mathbf{\mathbf{p_a}})$. Then we apply the rotation operator given by
$$
\mbox{ROTATE}(\mathbf{q_a}, \mathbf{r_{ab}}) = (\mathbf{r_{ab}})(\mathbf{q_a})(\mathbf{r_{ab}})^{-1} \, .
$$ 
The final position after the rotation is given by the \correct{imaginary part of the} unit quaternion $\mathbf{q_b} = (0, \mathbf{p_b})$. 

Introducing a parameterization $t$ such that $t=0$ and $t=1$ corresponding to the initial position $\mathbf{q_a}$ and $\mathbf{q_b}$, respectively, we can interpolate these two data points by the rotation operator $\mbox{ROTATE}(\mathbf{q_a}, \mathbf{r_{ab}}, t) = (\mathbf{r_{ab}})^t (\mathbf{q_a}) (\mathbf{r_{ab}})^{-t}$ for $t \in [0, 1]$. This expression leads to the so-called SLERP (\textit{Spherical Linear intERPolation}) formula \cite{shoemake_85,sola}: 
$$
\mbox{SLERP}(\mathbf{\mathbf{q_a}}, \mathbf{q_b}, t) = (\mathbf{q_a}) ((\mathbf{q_a})^{-1} \mathbf{q_b})^t, \quad t \in [0, 1] \, .
$$
In particular, if the quantity $\mathbf{a_{ab}}$ in the rotation quaternion is perpendicular (and this is assumed to be true in the remaining of this article) to both $\mathbf{\mathbf{p_a}}$ and $\mathbf{p_b}$, we have $\mathbf{\mathbf{p_a}} \cdot \mathbf{p_b} = \cos \theta_{ab}$ and 
\begin{equation*}
    \begin{split}
        \mathbf{p_b} &= (\cos \theta_{ab}) \mathbf{\mathbf{p_a}} + (\sin \theta_{ab}) (\mathbf{a_{ab}} \times \mathbf{\mathbf{p_a}}); \quad \quad \text{(Rodrigues' rotation formula \correct{\cite{rodrigues_40}})}\\
        \mathbf{a_{ab}} &= \mathbf{\mathbf{p_a}} \times \left[\dfrac{\mathbf{p_b} - (\cos \theta_{ab}) \mathbf{\mathbf{p_a}}}{(\sin \theta_{ab})}\right] = \dfrac{\mathbf{\mathbf{p_a}} \times \mathbf{p_b}}{\sin \theta_{ab}}.
    \end{split}
\end{equation*}

The function SLERP has two interesting properties. One is that the interpolant runs on the shortest path (geodesic) between both endpoints at a constant (angular) speed. If both two data points are both pure unit quaternions (i.e. $\mathbf{\mathbf{q_a}} = (0, \mathbf{\mathbf{p_a}})$ and $\mathbf{q_b} = (0, \mathbf{p_b})$), the SLERP interpolant between $\mathbf{q_a}$ and $\mathbf{q_b}$ lies on the vector part of the unit quaternion sphere. Mathematically, we can show that $\mbox{SLERP}(\mathbf{\mathbf{q_a}}, \mathbf{q_b}, t) = (0, \mathbf{p}(t))$ with $\lVert \mathbf{p}(t) \rVert = 1$ for all $t\in[0,1]$ if $\mathbf{p}(0) = \mathbf{\mathbf{p_a}}$ and $\mathbf{p}(1) = \mathbf{p_b}$. Another interesting property is that the inverse of a SLERP between two pure unit quaternions negates the interpolant, i.e. 
$$
(\mbox{SLERP}(\mathbf{\mathbf{q_a}}, \mathbf{q_b}, t))^{\pm 1} = \mbox{SLERP}(\pm \mathbf{\mathbf{q_a}}, \pm \mathbf{q_b}, t) = \pm \mbox{SLERP}(\mathbf{\mathbf{q_a}}, \mathbf{q_b}, t) \, .
$$

When we have more than two data points on a sphere, one might obtain a piecewise linear interpolant by applying SLERP in a piecewise fashion. If we want to obtain a smoother interpolant, one can consider the SQUAD (spherical quadrangle interpolation) \cite{dam_98}. It constructs a $C^1$ interpolant passing through the data points $\mathbf{p_i}$ and $\mathbf{p_{i + 1}}$ on $\mathbb{S}^2$ with two extra \correct{quaternions $\mathbf{s_i}$ and $\mathbf{s_{i+1}}$ as control points
\begin{eqnarray*}
\mathbf{s_{i}} &=& \mathbf{q_{i}} \exp\left[-\frac{1}{4} \left(\ln( \mathbf{q_{i}}^{-1} \mathbf{q_{i+1}}) + \ln(\mathbf{q_{i}}^{-1} \mathbf{q_{i-1}}) \right)\right] \\
\mathbf{s_{i+1}} &=& \mathbf{q_{i+1}} \exp\left[-\frac{1}{4} \left(\ln( \mathbf{q_{i+1}}^{-1} \mathbf{q_{i+2}}) + \ln(\mathbf{q_{i+1}}^{-1} \mathbf{q_{i}}) \right)\right]
\end{eqnarray*}
determined by $\mathbf{q_{i-1}}=(0,\mathbf{p_{i-1}})$, $\mathbf{q_{i}}=(0,\mathbf{p_{i}})$, $\mathbf{q_{i+1}}=(0,\mathbf{p_{i+1}})$, and $\mathbf{q_{i+2}}=(0,\mathbf{p_{i+2}})$.} In case $\mathbf{p_{i-1}}$ or $\mathbf{p_{i+2}}$ is not defined (in other words, we are at the beginning or the end of the data point sequence), they are defined as $\mathbf{p_i}$ or $\mathbf{p_{i+1}}$ respectively for convention. One first converts each $\mathbf{p_i}$ to a pure unit quaternion $\mathbf{q_i} = (0, \mathbf{\mathbf{p_i}})$, then the SQUAD interpolant is given by
\begin{eqnarray*}
& & \mbox{SQUAD}(\boxed{\underset{\text{for control}}{\mathbf{\mathbf{q_{i-1}}}}}, \boxed{\underset{\text{data}}{\mathbf{\mathbf{q_i}}}}, \boxed{\underset{\text{data}}{\mathbf{\mathbf{q_{i+1}}}}}, \boxed{\underset{\text{for control}}{\mathbf{\mathbf{q_{i+2}}}}}, t) \\
&=& \mbox{SLERP}(\mbox{SLERP}(\mathbf{\mathbf{q_i}}, \mathbf{q_{i + 1}}, t), \mbox{SLERP}(\mathbf{s_{i}}, \mathbf{s_{i + 1}}, t), 2t(1 - t)) \\
&=& \mathbf{q_i}\left(\mathbf{q_i}^{-1}\mathbf{q_{i+1}}\right)^t \left[\left[\mathbf{q_i}\left(\mathbf{q_i}^{-1}\mathbf{q_{i+1}}\right)^{t}\right]^{-1} \mathbf{s_{i}}\left( \mathbf{s_{i}}^{-1}\mathbf{s_{i + 1}}\right)^t\right]^{2t(1 - t)}  \correct{= \left(0,\mathbf{p}_{\mbox{SQUAD}}(t) \right)}
\end{eqnarray*}
for $t \in [0, 1]$. The interpolation SLERP and SQUAD are two most {widely} used interpolation schemes. There might be smoother explicit interpolations, but the expression could be complicated to implement and analyze. 


\section{Our Approaches}
\label{Sec:Proposal}

In this section, we will introduce two new ideas for interpolation on $\mathbb{S}^2$. In this first part of the discussion, we will introduce a simple interpolation method named \textit{Spherical Interpolation of orDER} $n$ (SIDER-$n-1$ in short) that gives a $C^{n-1}$ interpolant given $n \geq 3$ points on $\mathbb{S}^2$. The idea generalizes the construction of the B\'{e}zier curves developed for $\mathbb{R}$. However, similar to the standard polynomial interpolation, we would imagine that the interpolant will produce oscillations when the underlying function is not smooth enough. Therefore, instead of reconstructing the underlying interpolant using {single} highly smooth curve, we propose incorporating the ENO philosophy and developing a new \textit{Spherical Essentially Non-Oscillatory} (SENO) interpolation method. We will give the details in the second part of this section.

\subsection{SIDER Interpolations}

This section introduces a new class of interpolation schemes on the unit sphere, denoted by SIDER. We will first consider the case for $n=3$, i.e., three data points to obtain a $C^2$ curves on the unit sphere. We will then give a higher-order generalization for $n=4$ data points.

\begin{figure}
\centering
\includegraphics[width=0.4125\textwidth]{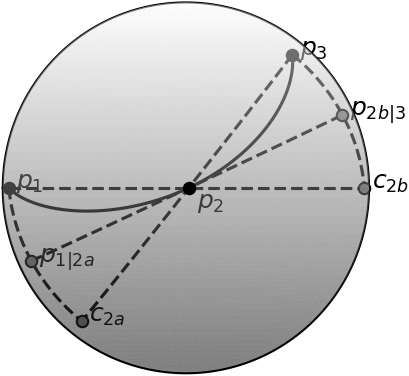}
\caption{Setup for SIDER2.}
\label{fig:sider2_ctrl_pts}
\end{figure}

With reference to the construction of quadratic B\'{e}zier curves, we propose the following spherical quadratic curve (denoted by SIDER2),
\begin{eqnarray*}
& &\mbox{SIDER2}(\boxed{\underset{\text{start}}{\mathbf{q_1}}}, \boxed{\underset{\text{second data}}{\mathbf{q_2}}}, \boxed{\underset{\text{end}}{\mathbf{q_3}}}, t) \\
&=&\mbox{SLERP}(\mbox{SLERP}(\mathbf{q_1},\mathbf{\mathbf{d_{2a}}}, t), \mbox{SLERP}(\mathbf{\mathbf{d_{2b}}}, \mathbf{q_3},t), t)\\
&=&\mathbf{q_1} \left(\mathbf{q_1}^{-1} {\mathbf{d_{2a}}}\right)^t \left[ \left[\mathbf{q_1}(\mathbf{q_1}^{-1} \mathbf{d_{2a}})^{t} \right]^{-1} \left[\mathbf{d_{2b}}(\mathbf{d_{2b}}^{-1} \mathbf{q_3})\right]^t\right]^{f_2(t)} \correct{= 
\left(0,\mathbf{p}_{\mbox{SIDER2}}(t) \right)}
\end{eqnarray*}
where $t \in [t_1, t_3]$, and $f_2(t) = (t - t_1)/(t_3 - t_1)$. Unless specified otherwise, we might use $t_1 = 0$ and $t_3 = 1$. The points $\mathbf{q_i} = (0, \mathbf{\mathbf{p_i}})$, $\mathbf{d_{2a}} = (0, \mathbf{\mathbf{c_{2a}}})$ and $\mathbf{d_{2b}} = (0, \mathbf{\mathbf{c_{2b}}})$ are the quaternion representation of the position vectors $\mathbf{\mathbf{p_i}}$, $\mathbf{\mathbf{c_{2a}}}$ and $\mathbf{\mathbf{c_{2b}}}$, respectively. We have shown the setup for the construction in Figure \ref{fig:sider2_ctrl_pts}.

The crucial step in this expression is how we determine the control points $\mathbf{c_{2a}}$ and $\mathbf{c_{2b}}$. We construct $\mathbf{c_{2b}}$ (and $\mathbf{c_{2a}}$) using the geodesic extrapolating based on the first data points $\mathbf{p_1}$ (and $\mathbf{p_3}$) and the intermediate one $\mathbf{p_2}$ so that the final interpolant reaches $\mathbf{p_2}$ when $t = t_2 = 0.5(t_1 + t_3)$. To enforce this condition, we refer to the spatial relationships among the data points and the (only) control point in a quadratic force interpolating B\'{e}zier curve. Mathematically, we assign 
$$
\mathbf{d_{2a}} = (0, \mathbf{\mathbf{c_{2a}}}) = \text{SLERP}(\mathbf{q_3},\mathbf{q_2},2) \, \mbox{ and } \,
\mathbf{d_{2b}} = (0, \mathbf{\mathbf{c_{2b}}}) = \text{SLERP}(\mathbf{q_1}, \mathbf{q_2},2) \, .
$$
Although SIDER2 seems to directly replace the linear interpolation in the quadratic B\'{e}zier curve in ${\mathbb{R}^n}$ by the SLERP, these two expressions are different. A typical B\'{e}zier interpolant usually passes through only the first and the last points while treating all other points as control points. On the other hand, our SIDER2 is an interpolation formula that requires the interpolant to pass through all given sampling points.

Here we give several properties of SIDER2. 
\begin{itemize}
\item The interpolant passes through all three data points. We can easily check that 
\[\mbox{SIDER2}(\mathbf{q_1},\mathbf{q_2},\mathbf{q_3},0)=(0, \mathbf{p_1}) = \mathbf{q_1} \mbox{ and } \mbox{SIDER2}(\mathbf{q_1},\mathbf{q_2},\mathbf{q_3},1)=(0, \mathbf{p_3}) = \mathbf{q_3} \, .\]
For the second data point, we have 
\begin{eqnarray*}
\text{SIDER2}(\mathbf{q_1}, \mathbf{q_2},\mathbf{q_3},0.5) &=& \text{SLERP}(\text{SLERP}(\mathbf{q_1}, \mathbf{\mathbf{d_{2a}}}, 0.5), \text{SLERP}(\mathbf{\mathbf{d_{2b}}}, \mathbf{q_3},0.5), 0.5) \\
&=& \text{SLERP}(\mathbf{q_{1|2a}}, \mathbf{q_{2b|3}}, 0.5) \, ,
\end{eqnarray*} where \correct{$\mathbf{q_{1|2a}}=(0,\mathbf{p_{1|2a}})$, $\mathbf{q_{2b|3}}=(0,\mathbf{p_{2b|3}})$,} $\mathbf{p_{1|2a}}$ is the midpoint along the geodesic between $\mathbf{p_1}$ and $\mathbf{c_{2a}}$, and $\mathbf{p_{2b|3}}$ is the midpoint along the geodesic between $\mathbf{c_{2b}}$ and $\mathbf{p_3}$. The last expression looks for the midpoint along the geodesic between $\mathbf{p_{1|2a}}$ and $\mathbf{p_{2b|3}}$ which, therefore, leads to $\mathbf{p_2}$ according to the construction of $\mathbf{c_{2a}}$ and $\mathbf{c_{2b}}$.
\item Since SLERP guarantees that the interpolant stays on the sphere, we are assured that SIDER2 also generates points with the unit length for all $t\in[t_1,t_3]$. 
\item Reversing the start and end points (while keeping the midway data obtained at $t_2 = 0.5(t_1 + t_3)$) would not change the curve on the sphere.
\end{itemize}

The approach can be easily extended to the case where $n=4$, i.e., we construct a $C^3$ curve from four given data points. The construction is given by
\begin{eqnarray*}
&& \mbox{SIDER3}(\boxed{\underset{\text{start}}{\mathbf{q_1}}}, \boxed{\underset{\text{second point}}{\mathbf{q_2}}}, \boxed{\underset{\text{third point}}{\mathbf{q_3}}}, \boxed{\underset{\text{end}}{\mathbf{q_4}}}, t) \\
&=& \mbox{SLERP}(\mbox{SIDER2}(\mathbf{q_1}, \mathbf{q_2},\mathbf{q_3},g_3(t)), \mbox{SIDER2}(\mathbf{q_2},\mathbf{q_3},\mathbf{q_4},h_3(t)), f_3(t)) \correct{= 
\left(0,\mathbf{p}_{\mbox{SIDER3}}(t) \right)}
\end{eqnarray*}
where $t \in [t_1, t_4]$. Therefore, a SIDER3 reconstruction is a linear combination of two scaled SIDER2, that we interpolate within $\{\mathbf{p_1}, \mathbf{p_2}, \mathbf{p_3}\}$ and $\{\mathbf{p_2}, \mathbf{p_3}, \mathbf{p_4}\}$ simultaneously.

Given that the timestamps are equalized, i.e. $t_2 = \frac{1}{3}(2t_1 + t_4)$ and $t_3 = \frac{1}{3}(t_1 + 2t_4)$, we want the functions $g_3(t)$, $h_3(t)$ and $f_3(t)$ satisfy the conditions imposed at $t_1$, $t_2$, $\frac{1}{2}(t_1+t_4)$, $t_3$ and $t_4$ as shown in Table \ref{Table:SIDER3}. The simplest linear and quadratic functions that satisfy these constraints are given by $g_3(t) = (t - t_1)/(t_3 - t_1)$, $h_3(t) = (t - t_2)/(t_4 - t_2)$ and $f_3(t) = (t - t_1)/(t_4 - t_1)$. For simplicity, we set the starting time $t_1 = 0$ and the ending time $t_4 = 1$, so that $t_2 = \frac{1}{3}$ and $t_3 = \frac{2}{3}$, and $g_3(t) = 3t/2$, $h_3(t) = (3t - 1)/2$ and $f_3(t) = t$. 

\begin{table}
\centering
\begin{tabular}{|c||c|c|c|}
\hline 
$t$                    & $g_3(t)$      & $h_3(t)$      & $f_3(t)$      \\ 
\hline\hline 
$t_1$                  & 0             &              & 0             \\ 
\hline 
$t_2 := \frac{1}{3}(2t_1 + t_4)$                  & 1/2 & 0 & 1/3 \\ 
\hline 
$(t_1+t_4)/2$ &  &  & 1/2 \\ 
\hline 
$t_3 := \frac{1}{3}(t_1 + 2t_4)$                  & 1 & 1/2 & 2/3 \\ 
\hline 
$t_4$                  &               & 1             & 1      \\      
\hline 
\end{tabular}
\caption{Conditions imposed on the functions $g_3(t)$, $h_3(t)$ and $f_3(t)$ in the construction of SIDER3. We leave the box blank if we do not {account for the function value} at that timestamp.} 
\label{Table:SIDER3}
\end{table}

It is easy to check that the interpolant satisfies
$$
\mbox{SIDER3} \left(\mathbf{q_1},\mathbf{q_2},\mathbf{q_3},\mathbf{q_4},0\right)={(0, \mathbf{p_1}) = \mathbf{q_1}} \, \mbox{ and } \,
\mbox{SIDER3} \left(\mathbf{q_1},\mathbf{q_2},\mathbf{q_3},\mathbf{q_4},1\right)={(0, \mathbf{p_4}) = \mathbf{q_4}} \, .
$$
For the two other intermediate data points, we have 
\begin{eqnarray*}
& &\mbox{SIDER3} \left(\mathbf{q_1},\mathbf{q_2},\mathbf{q_3},\mathbf{q_4},\frac{1}{3} \right) \\
&=& \mbox{SLERP}\left(\mbox{SIDER2}\left(\mathbf{q_1},\mathbf{q_2},\mathbf{q_3},g_3\left(\frac{1}{3}\right)\right), \mbox{SIDER2}\left(\mathbf{q_2},\mathbf{q_3},\mathbf{q_4},h_3\left(\frac{1}{3}\right)\right), f_3\left(\frac{1}{3}\right)\right) \\
&=& \mbox{SLERP}\left(\mbox{SIDER2}\left(\mathbf{q_1},\mathbf{q_2},\mathbf{q_3},1/2\right), \mbox{SIDER2}\left(\mathbf{q_2},\mathbf{q_3},\mathbf{q_4},0\right), 1/3\right) \\
&=& {\mbox{SLERP}\left(\mathbf{q_2},\mathbf{q_2},\frac{1}{3}\right) = {(0, \mathbf{p_2}) = \mathbf{q_2}}} \, . 
\end{eqnarray*} 
and
\begin{eqnarray*}
& &\mbox{SIDER3} \left(\mathbf{q_1},\mathbf{q_2},\mathbf{q_3},\mathbf{q_4},\frac{2}{3} \right) \\
&=& \mbox{SLERP}\left(\mbox{SIDER2}\left(\mathbf{q_1},\mathbf{q_2},\mathbf{q_3},g_3\left(\frac{2}{3}\right)\right), \mbox{SIDER2}\left(\mathbf{q_2},\mathbf{q_3},\mathbf{q_4},h_3\left(\frac{2}{3}\right)\right), f_3\left(\frac{2}{3}\right)\right) \\
&=& \mbox{SLERP}\left(\mbox{SIDER2}\left(\mathbf{q_1},\mathbf{q_2},\mathbf{q_3},1\right), \mbox{SIDER2}\left(\mathbf{q_2},\mathbf{q_3},\mathbf{q_4},1/2\right), 2/3\right) \\
&=& {\mbox{SLERP}\left(\mathbf{q_3},\mathbf{q_3},\frac{2}{3}\right) = {(0, \mathbf{p_3}) = \mathbf{q_3}}} \, . 
\end{eqnarray*}
Even though it is possible to derive the time derivatives of SIDER3 analytically, the expressions are complex, and we do not present them explicitly here. However, in the example section, we will numerically demonstrate that the interpolant and {its first few} derivatives are all continuous.

It is also possible to further develop higher-order SIDER interpolants using recursions. In general, we can construct a SIDER-$n$ interpolant to interpolate $n+1$ data points with equal time spacing from $t = 0$ to $t = 1$. The expression consists of one SLERP of two SIDER-$(n-1)$ interpolants with $g_n(t) = nt/(n-1), h_n(t) = g_n(t) - 1/(n-1)$ and $f_n(t) = t$. For example, we might construct a SIDER4 formula based on one SLERP and two SIDER3's,
\begin{eqnarray*}
& & \text{SIDER4}(\mathbf{q_1},\mathbf{q_2},\mathbf{q_3},\mathbf{q_4},\mathbf{q_5}, t) \\
&=& \text{SLERP}\left(\text{SIDER3}\left(\mathbf{q_1},\mathbf{q_2},\mathbf{q_3},\mathbf{q_4},g_4(t)\right), \text{SIDER3}\left(\mathbf{q_2},\mathbf{q_3},\mathbf{q_4},\mathbf{q_5}, h_4(t)\right), f_4(t)\right) \, ,
\end{eqnarray*} where $g_4(t) = {(4/3)t}, h_4(t) = g_4(t) - 1/3$, $f_4(t) = {t}$ with equal time spacing data.

\subsection{Constraints on the Data Points}

Unlike any typical B\'{e}zier construction in the Euclidean space, the above approach might not work for some datasets. This section will discuss some constraints on the given dataset for our SIDER {reconstructions.} 

In the quaternion representation of points on the sphere, opposite points (i.e., $\mathbf{p_1}$ and $-\mathbf{p_1}$) are equivalent and are treated as the same rotation operator. To avoid the non-uniqueness in our algorithms, we need to constrain the sector angle between any two adjacent data points, denoted by $\theta$, to be less than $\pi/2$. This condition also implies that the geodesic distance between adjacent data points must be smaller than $\pi/2$. Otherwise, our algorithm will move the data points to their opposite side so that the geodesic distance becomes smaller than $\pi/2$. If the geodesic distance between {adjacent points $\mathbf{p_i}$ and $\mathbf{p_{i+1}}$} is exactly $\pi/2$, the same is true for the geodesic distance between {$-\mathbf{p_i}$ and $\mathbf{p_{i+1}}$. Therefore, we cannot uniquely decide whether $\mathbf{p_i}$ or $-\mathbf{p_i}$ is closer to $\mathbf{p_{i+1}}$, so we cannot give a unique interpolant from $\mathbf{p_i}$ to $\mathbf{p_{i+1}}$. This leads to the location of the deduced control points undetermined, and the interpolant might not look like what the users would expect.}

Another constraint is that all the data points should not lie on the same great circle. Otherwise, any high order interpolation on the sphere will be degenerated, as in the Euclidean space where we have colinear data points.

\subsection{Spherical-ENO (SENO) Interpolation Based on SIDER}

One usually observes oscillations in the interpolant when reconstructing a high-order curve with sharp changes and turns, and this behavior is undesirable in many applications. This section follows the philosophy of Essentially Non-Oscillatory (ENO) and proposes an ENO interpolation on the unit sphere. We name the interpolation approach the \textit{Spherical Essentially Non-Oscillatory} (SENO in short).

Given a set of $2n$ data points, denoted by $\mathbf{p_{i-n+1}}, \cdots, \mathbf{p_i}$, $\mathbf{p_{i+1}}, \cdots \mathbf{p_{i+n-1}}$, we are interested in constructing a high-order curve between $\mathbf{p_i}$ and $\mathbf{p_{i+1}}$. To do this, we first reconstruct a $C^{n}$ curve from any $n+1$ consecutive data points using SIDER-$n$. For example, for $n=2$, i.e., we are given four data points, we first construct two $C^2$ SIDER2 curves from any three consecutive data on the unit sphere. When $n=3$, we have in total six data points. From these, we obtain three $C^3$ curves obtained by SIDER3. To avoid an oscillatory interpolant, we consider these $n$ interpolants from {SIDER}-$(n-1)$ and determine the corresponding \textit{variation} of these curves between the data points $\mathbf{p_{i}}$ and $\mathbf{p_{i+1}}$. The one with the least \textit{variation} is chosen to represent the SENO interpolant between the points $\mathbf{p_i}$ and $\mathbf{p_{i+1}}$.

There are various ways how one can define the \textit{variation} of a curve on the unit sphere. For simplicity, we define it to be the length of the segment joining $\mathbf{p_i}$ and $\mathbf{p_{i+1}}$. With the minimum length achieved by the geodesic, a longer length indicates a larger variation and more violent oscillations obtained in the high-order reconstruction. Instead of determining an analytical expression for the variation by evaluating the exact line integral, we propose applying the composite Trapezoidal rule to numerically approximate the distance between these data points. Between the points $\mathbf{p_i}$ and $\mathbf{p_{i+1}}$, we insert $k$ intermediate points ($k=3$ in all numerical simulations we have obtained) so that we have $k+1$ small segments of equal time span. Then we approximate the total length of the curve by summing up the geodesic distance of these $k+1$ small segments. {This} underestimates the total length of the interpolant from $\mathbf{p_i}$ and $\mathbf{p_{i+1}}$ {but} we observe that such approximation already provides a reasonable choice of SIDER-$n$ for a non-oscillatory reconstruction.

We now consider the computational efficiency. Both SQUAD and SIDER2 interpolations require three SLERP reconstructions. However, since SENO2 constructs two SIDER2 interpolants, the overall computational complexity for a SENO2 is double {to} that of a SQUAD. Since a SIDER3 reconstruction requires one call of SLERP and two calls of SIDER2, it requires seven SLERP calls in total. For SENO3, there are three different sets of stencils producing three SIDER3 interpolants. Therefore, the total number of SLERP calls is 21 for a single SENO3 interpolation, implying seven times the operations of a SQUAD interpolation. 

To end the discussion, we summarize several properties of both SIDER and SENO in Table \ref{Table:Comparison}.

\begin{table}[!ht]
\centering
\begin{tabular}{|c||cc|cc|cc|}
\hline
& SLERP & SQUAD & SIDER2 & SENO2 & SIDER3 & SENO3 \\
\hline
\hline
Data points & 2 & 4 & 3 & 4 & 4 & 6 \\
\hline
Order (for smooth curves) & 2 & 3 & 3 & 3 & 4 & 4 \\ \hline
Complexity (SLERP calls) & 1 & 3 & 3 & 6 & 7 & 21 \\ \hline
\end{tabular}
\caption{Different properties of SLERP, SQUAD, SIDER2, SENO2, SIDER3 and SENO3 including the total number of data points required, the expected convergence order if the underlying curve is smooth enough, and the computational complexity in terms of the SLERP calls.}
\label{Table:Comparison}
\end{table}


\section{Numerical Examples}
\label{Sec:Example}

This section will consider several examples to compare our proposed SIDER and SENO interpolations. We will compare our reconstructions with those from (piecewise) SLERP and SQUAD to show that our approach can achieve a high-order reconstruction while avoiding over-shooting when the underlying curve contains some sharp turns. 


\begin{figure}[!ht]
\centering
(a)\includegraphics[width=0.425\textwidth]{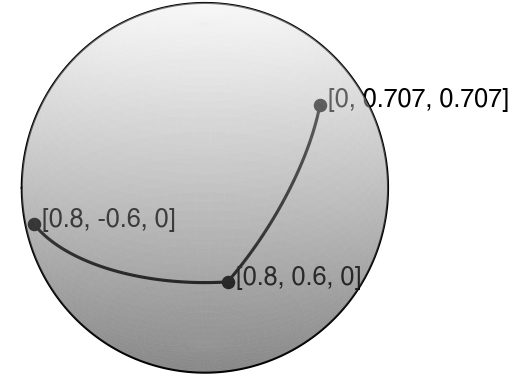}
(b)\includegraphics[width=0.425\textwidth]{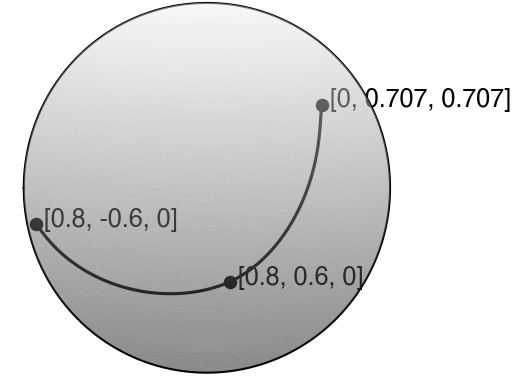} \\
(c)\includegraphics[width=0.425\textwidth]{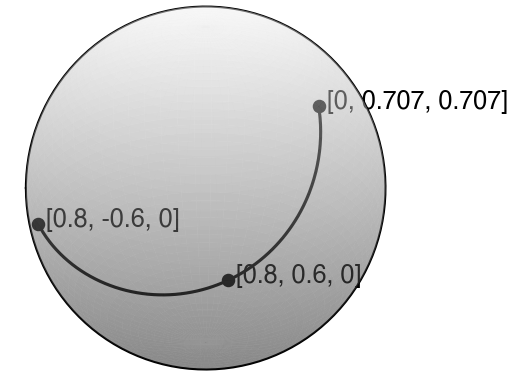}
(d)\includegraphics[width=0.425\textwidth]{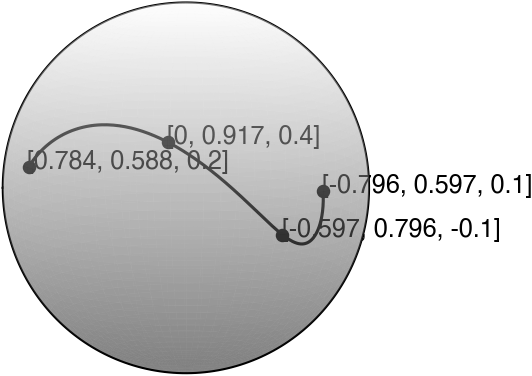}
\caption{(Section \ref{Ex:SimpleCase}) Interpolation of three data points based on (a) SLERP, (b) SQUAD, (c) SIDER2, {and of four data points based on (d) SIDER3.}}
\label{Fig:SimpleCase}
\end{figure}

\begin{figure}[!ht]
\centering
(a)\includegraphics[width=0.925\textwidth,clip]{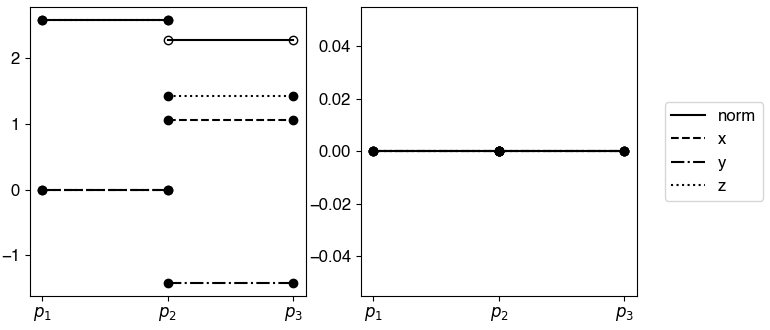}
(b)\includegraphics[width=0.925\textwidth,clip]{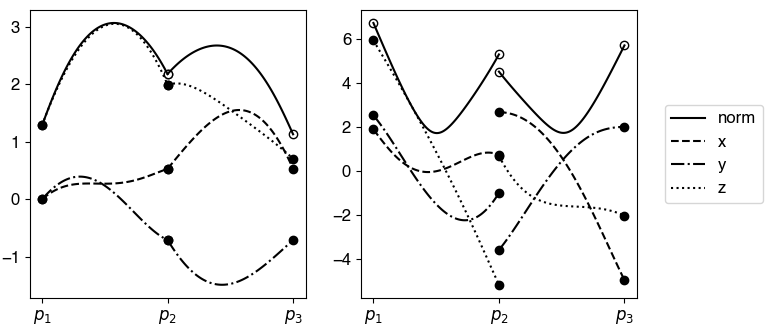}
(c)\includegraphics[width=0.925\textwidth,clip]{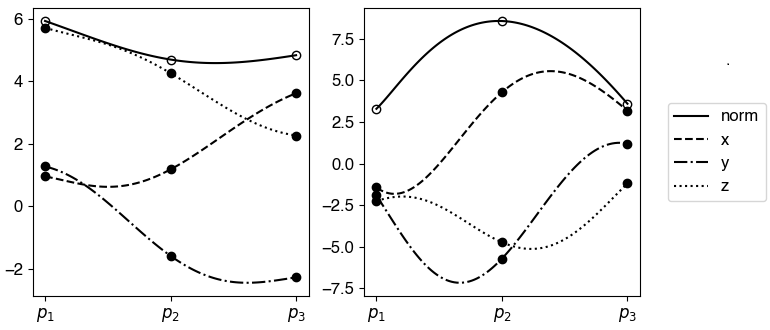}
\caption{(Section \ref{Ex:SimpleCase}) The (left) angular velocity and the (right) acceleration computed from the interpolant given by (a) SLERP, \correct{(b) SQUAD}, and (c) SIDER2. Since the derivatives are also quaternions, we plot the norm, and the individual components of the vector part of the corresponding quaternion separately in each subplot. Since the real part of all derivatives is zero, they are all omitted here.}
\label{Fig:SimpleCase_Component}
\end{figure}

\begin{table}[!ht]
\centering
(a) 
    \begin{tabular}{|c||c|c|}
        \hline
        point & location & $t$ \\
        \hline\hline
        $\mathbf{p_1}$ & $(0.8, -0.6, 0)$ & 0 \\
        \hline
        $\mathbf{p_2}$ & $(0.8, 0.6, 0)$ & 0.5 \\
        \hline
        $\mathbf{p_3}$ & $(0, \sqrt{0.5}, \sqrt{0.5})$ & 1 \\
        \hline
    \end{tabular} \hspace{0.5cm}
(b)
    \begin{tabular}{|c||c|c|}
        \hline
        point & location & $t$ \\
        \hline\hline
        $\mathbf{p_1}$ & $(\sqrt{0.6144}, \sqrt{0.3456}, 0.2)$ & 0 \\
        \hline
        $\mathbf{p_2}$ & $(0, \sqrt{0.84}, 0.4)$ & 1/3 \\
        \hline
        $\mathbf{p_3}$ & $(-\sqrt{0.3564}, \sqrt{0.6336}, -0.1)$ & 2/3 \\
        \hline
        $\mathbf{p_4}$ & $(-0.64, 0.48, 0.6)$ & 1 \\
        \hline
    \end{tabular}
\caption{(Section \ref{Ex:SimpleCase}) (a) Data points for SLERP, SQUAD and SIDER2. (b) Data points for SIDER3.}
\label{tab:test_case_14}
\end{table} 

\begin{figure}[!ht]
\centering
\includegraphics[width=0.99\textwidth]{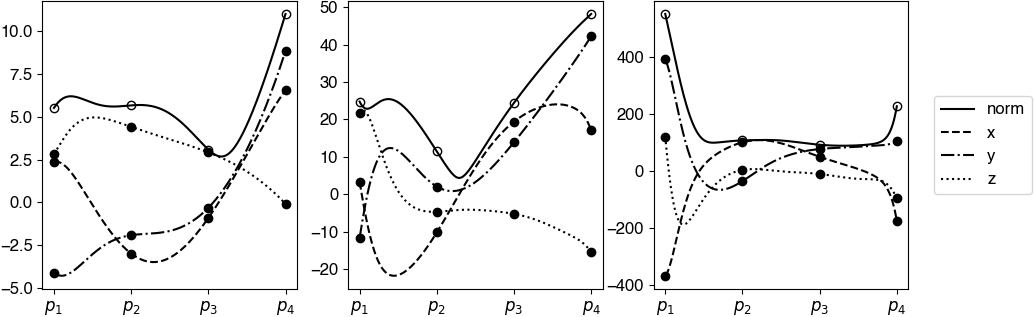}
\caption{(Section \ref{Ex:SimpleCase}) (From left to right) The first three derivatives of the interpolant computed from SIDER3. Since the derivatives are also quaternions, we plot the norm, and the individual components of the vector part of the corresponding quaternion separately in each subplot. Since the real part of all derivatives is zero, they are all omitted here.} 
\label{fig:sider3_eg4_ddt}
\end{figure}

\subsection{Comparison of SLERP, SQUAD, and SIDER}
\label{Ex:SimpleCase}

In this simple example, we interpolate the data points given in Table \ref{tab:test_case_14}(a) using SLERP, SQUAD, and SIDER2. For the SLERP reconstruction, we apply the SLERP to any adjacent data point. This process essentially gives a piecewise SLERP reconstruction. For SQUAD, we also interpolate the data points in a piecewise fashion. The control points of the SQUAD would follow the so-called \textit{bilinear parabolic blending} \cite{haarbach_18}, so most intermediate data point would appear in four consecutive SQUAD pieces, i.e. given a sequence $\{\mathbf{p_1}, \mathbf{p_2}, \mathbf{p_3}, \cdots, \mathbf{p_n}\}$, every data point $\mathbf{p_i}$, where $i = 3$ to $n-2$ will appear in the $(i-2)^{\text{th}}$, $(i-1)^{\text{th}}$, $i^{\text{th}}$ and $(i+1)^{\text{th}}$ interpolants. As for using SQUAD, we let $\mathbf{p_0} = \mathbf{p_1}$ and $\mathbf{p_4} = \mathbf{p_3}$ as we calculate the control points unless extra data points are given like in the simulation performed in Section \ref{Sec:AccuracyConvergence}.

Figure \ref{Fig:SimpleCase}(a-c) plot the interpolation results on the sphere based on SLERP, SQUAD and SIDER2 interpolations. Figure \ref{Fig:SimpleCase_Component} shows the corresponding angular derivatives of all interpolants as a function of $t$. Since the piecewise SLERP interpolates any two adjacent data points along the geodesics, the interpolant has sharp kinks at the data point $\mathbf{p_2}$. The first derivative of the interpolant (i.e., the angular velocity) is piecewise constant, as shown in the left subplot of Figure \ref{Fig:SimpleCase_Component}(a). The corresponding higher derivatives (i.e., the {angular} acceleration and others) are zero. Even though SQUAD can determine a much smoother interpolant than the SLERP, the reconstruction provides only a $C^1$ curve. As shown on the right subplot in Figure \ref{Fig:SimpleCase_Component}(b), we can see that the acceleration is not continuous at the middle data point. This property indicates that the interpolant cannot achieve a higher regularity than $C^1$. Our proposed SIDER2, on the other hand, not only shows a smooth reconstruction in Figure \ref{Fig:SimpleCase}(c), both the angular velocity and the acceleration are continuous, \correct{as shown in Figure \ref{Fig:SimpleCase_Component}(c)}. This observation shows that our reconstruction can produce a $C^2$ curve. The analytical expressions of the (angular) derivatives of SQUAD and SIDER2 are tedious. We give the derivations in the appendix for reference.

With one more point on the sphere, we can interpolate the data using SIDER3. We consider the dataset given in Table \ref{tab:test_case_14}(b). The resulting interpolant is shown in Figure \ref{Fig:SimpleCase}(d). We observe that the first three derivatives of the interpolant are all continuous, as shown in Figure \ref{fig:sider3_eg4_ddt}. 


\begin{table}[!ht]
    \centering
    \begin{tabular}{|c|c|c|}
        \hline
        Test case & (a) & (b) \\
        \hline
        \hline
        $\mathbf{p_1}$ & \multicolumn{2}{c|}{$(\sqrt{0.6144}, \sqrt{0.3456}, 0.2)$} \\
        \hline
        $\mathbf{p_2}$ & \multicolumn{2}{c|}{$(0, \sqrt{0.84}, 0.4)$} \\
        \hline
        $\mathbf{p_3}$ & \multicolumn{2}{c|}{$(-\sqrt{0.3564}, \sqrt{0.6336}, -0.1)$} \\
        \hline
        $\mathbf{p_4}$ & $(-0.64, 0.48, 0.6)$ & $(-\sqrt{0.6336}, \sqrt{0.3564}, 0.1)$ \\
        \hline
    \end{tabular}
    \caption{(Section \ref{SubSec:Ex23}) Two sets of data to demonstrate the effectiveness of SENO2.}
    \label{tab:case2and3}
\end{table}

\begin{figure*}[!ht]
\centering
(a)\includegraphics[width=0.45\textwidth]{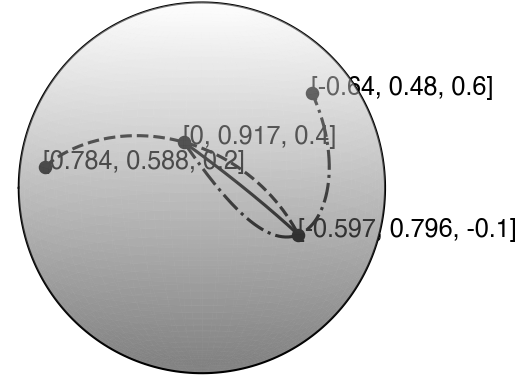}
(b)\includegraphics[width=0.45\textwidth]{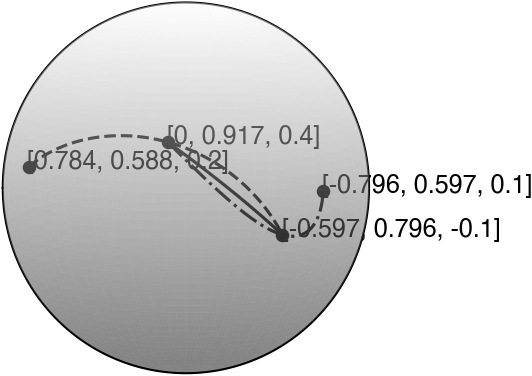}
\caption{(Section \ref{SubSec:Ex23}) The interpolant by SENO2 of the dataset (a) and (b) as defined in Table \ref{tab:case2and3}. Both candidates $S_{123}$ (the dashed line) and $S_{234}$ (the dashed-dotted line) are constructed from SIDER2. The SLERP between $\mathbf{p_2}$ and $\mathbf{p_3}$ is also plotted (in solid line) for reference.}
\label{fig:Ex2}
\end{figure*}

\subsection{Reconstruction by SENO2}
\label{SubSec:Ex23}
This example demonstrates the idea in SENO2 constructed based on SIDER2. Table \ref{tab:case2and3} shows the two sets of data that we consider. In both examples, we see that the last data point creates a sharp turn at the point $\mathbf{p_3}$. When constructing a highly smooth reconstruction, we expect the interpolant between $\mathbf{p_2}$ and $\mathbf{p_3}$ to oscillate and be affected by the kink at $\mathbf{p_3}$. In each of the data sets, we, therefore, construct two SIDER2 interpolants $S_{123}$ and $S_{234}$ from four given data points. We show the solutions from our SIDER2 in Figure \ref{fig:Ex2}. In both subplots, we have also shown the piecewise SLERP interpolant between $\mathbf{p_2}$ and $\mathbf{p_3}$ using solid lines. From these two tests, we see that the SENO2 reconstruction from the test case (a) in Table \ref{tab:case2and3} would choose $S_{123}$ while test case (b) would choose $S_{234}$ because their {\textit{variation}} is smaller than the other one between the points in concern, i.e. $\mathbf{p_2}$ and $\mathbf{p_3}$.


\subsection{Reconstruction by SENO3}
\label{SubSec:Ex5}

\begin{table}[!ht]
\centering
\begin{tabular}{|c|c|}
\hline
$\mathbf{p_1}$ & $(-0.9462408024134863, 0.2340693569139826, -0.2232484714432692)$ \\
\hline
$\mathbf{p_2}$ & $(-0.5756591575040059, 0.7203584217199284, -0.3869112025244969)$ \\
\hline
$\mathbf{p_3}$ & $(-0.5139135508439371, 0.8072140040848369, 0.29034189134243293)$ \\
\hline
$\mathbf{p_4}$ & $(0.1733822829796129, 0.5285757390277231, 0.830991138376381)$ \\
\hline
$\mathbf{p_5}$ & $(0.8196895318805648, -0.045366259610012546, 0.571008733571053)$ \\
\hline
$\mathbf{p_6}$ & $(0.8410803457569805, 0.5409102069487302, 0)$ \\
\hline
\end{tabular}
\caption{(Section \ref{SubSec:Ex5}) Data to demonstrate the effectiveness of SENO3.}
\label{tab:case5and6}
\end{table}

\begin{figure}[!ht]
\centering
\includegraphics[width=0.5\textwidth]{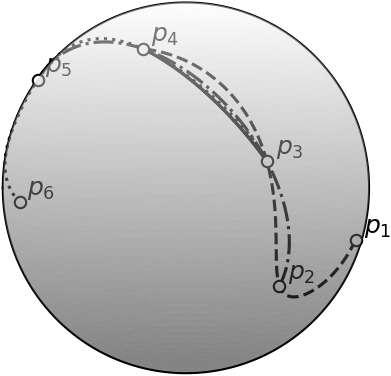}
\caption{(Section \ref{SubSec:Ex5}) The interpolant by SENO3. All three candidates $S_{1234}$ (the dashed line), $S_{2345}$ (the dashed-dotted line) and $S_{3456}$ (the dotted line) are constructed from SIDER3. The SLERP between $\mathbf{p_3}$ and $\mathbf{p_4}$ is also plotted (in solid line) for reference.}
\label{fig:sbzp3_select}
\end{figure}

In this example, we demonstrate the effectiveness of SENO3 constructed based on SIDER3. Table \ref{tab:case5and6} shows the dataset containing six \correct{random} data points on the sphere. We apply SIDER3 to reconstruct a $C^3$ curve for any four consecutive data points. This step leads to three different SIDER3 interpolants, denoted by $S_{1234}, S_{2345}$ and $S_{3456}$. For this case with six data points, SENO3 is interested in obtaining an interpolant that {varies the least} between $\mathbf{p_3}$ and $\mathbf{p_4}$. This condition implies the interpolant closest to the geodesic (as obtained from SLERP) joining $\mathbf{p_3}$ to $\mathbf{p_4}$. Figure \ref{fig:sbzp3_select} shows all SIDER3 interpolants together with the SLERP from $\mathbf{p_3}$ to $\mathbf{p_4}$ (in solid lines). As we can see from this figure, SENO3 will pick the dotted line corresponding to the SIDER3 connecting $\mathbf{p_3}$ to $\mathbf{p_6}$, to represent the high-order interpolant between $\mathbf{p_3}$ and $\mathbf{p_4}$.

\begin{table}[!htb]
\scriptsize
\centering
(a)
\begin{tabular}{|c||cc|cc|cc|cc|}
\hline
$1/\Delta t$ & $e_{\text{SLERP}}$ & $\rho_{\Delta t, \text{SLERP}}$ & $e_{\text{SQUAD}}$ & $\rho_{\Delta t, \text{SQUAD}}$ & $e_{\text{SENO2}}$ & $\rho_{\Delta t, \text{SENO2}}$ & $e_{\text{SENO3}}$ & $\rho_{\Delta t, \text{SENO3}}$ \\
\hline
\hline
16           & 7.5383e-03         & -                               & 5.3563e-03         & -                               & 3.3255e-03         & -                               & 5.7331e-03         & -                               \\
32           & 1.9053e-03         & 1.9842                          & 1.1249e-03         & \textbf{2.2514}                          & 6.0749e-04         & 2.4526                          & 2.2092e-04         & 4.6977                          \\
64           & 4.7560e-04         & 2.0022                          & 2.6434e-04         & \textbf{2.0894}                          & 7.6428e-05         & 2.9907                          & 1.2270e-05         & 4.1703                          \\
128          & 1.1909e-04         & 1.9977                          & 6.4493e-05         & \textbf{2.0352}                          & 9.2140e-06         & 3.0522                          & 6.9148e-07         & 4.1493                          \\
256          & 2.9761e-05         & 2.0006                          & 1.5909e-05         & \textbf{2.0193}                          & 1.1237e-06         & 3.0356                          & 4.1101e-08         & 4.0724                          \\
512          & 7.4348e-06         & 2.0011                          & 3.9474e-06         & \textbf{2.0108}                         & 1.3888e-07         & 3.0164                          & 2.5162e-09         & 4.0299                          \\
1024         & 1.8532e-06         & 2.0042                          & 9.8288e-07         & \textbf{2.0058}                         & 1.7248e-08         & 3.0093                          & 1.5658e-10         & 4.0063                          \\
2048         & 4.5786e-07         & 2.0171                          & 2.4521e-07         & \textbf{2.0030}                         & 2.1491e-09         & 3.0047                          & 9.7544e-12         & 4.0047                          \\
4096         & 1.0901e-07         & 2.0704                          & 6.1236e-08         & \textbf{2.0015}                          & 2.6820e-10         & 3.0024                          & 6.0906e-13         & 4.0014                         
\\
\hline
\end{tabular}
\medskip \\
(b)
\begin{tabular}{|c||cc|cc|cc|cc|}
\hline
$1/\Delta t$ & $e_{\text{SLERP}}$ & $\rho_{\Delta t, \text{SLERP}}$ & $e_{\text{SQUAD}}$ & $\rho_{\Delta t, \text{SQUAD}}$ & $e_{\text{SENO2}}$ & $\rho_{\Delta t, \text{SENO2}}$ & $e_{\text{SENO3}}$ & $\rho_{\Delta t, \text{SENO3}}$ \\
\hline
\hline
16           & 7.5383e-03         & -                               & 2.2475e-03         & -                               & 5.3375e-03         & -                               & 2.6273e-03         & -                               \\
32           & 1.9053e-03         & 1.9842                          & 2.0176e-04         & \textbf{3.4776}                          & 6.6980e-04         & 2.9944                          & 1.7398e-04         & 3.9166                          \\
64           & 4.7560e-04         & 2.0022                          & 2.0375e-05         & \textbf{3.3078}                          & 7.7677e-05         & 3.1082                          & 1.0571e-05         & 4.0407                          \\
128          & 1.1909e-04         & 1.9977                          & 2.3083e-06         & \textbf{3.1419}                          & 9.2349e-06         & 3.0723                          & 6.5747e-07         & 4.0071                          \\
256          & 2.9761e-05         & 2.0006                          & 2.7846e-07         & \textbf{3.0513}                          & 1.1240e-06         & 3.0384                          & 4.0534e-08         & 4.0197                          \\
512          & 7.4348e-06         & 2.0011                          & 3.4414e-08         & \textbf{3.0164}                          & 1.3888e-07         & 3.0167                          & 2.5072e-09         & 4.0150                          \\
1024         & 1.8532e-06         & 2.0042                          & 4.2869e-09         & \textbf{3.0050}                          & 1.7248e-08         & 3.0094                          & 1.5644e-10         & 4.0024                          \\
2048         & 4.5786e-07         & 2.0171                          & 5.3526e-10         & \textbf{3.0016}                          & 2.1491e-09         & 3.0047                          & 9.7522e-12         & 4.0037                          \\
4096         & 1.0901e-07         & 2.0704                          & 6.6881e-11         & \textbf{3.0006}                          & 2.6820e-10         & 3.0024                          & 6.0903e-13         & 4.0012                         
\\
\hline
\end{tabular}
\caption{\correct{(Section \ref{Sec:AccuracyConvergence}) The error and the numerical accuracy $\rho_{\Delta t}$ for SIDER, SQUAD, SENO2 and SENO3 when the generating function is (a) non-differentiable and (b) smooth.}} 
\label{tab:area_err}
\end{table}

\subsection{Accuracy and Convergence}
\label{Sec:AccuracyConvergence}

This section performs some numerical tests to verify the numerical accuracy of the proposed interpolation schemes. We will check the convergence of the methods as we insert more data points for curve reconstruction. We first construct the following {two} curves with different smoothness. \correct{
For the smooth curve, we first construct a parametrized curve given by $y(t)=t$ and 
$$
z(t)=f(t)=\exp\left( - \frac{t^2}{2\sigma^2} \right) \sin (2\pi t)
$$
with $\sigma=0.1$ on the $x=1$ plane. This definition gives an infinitely differentiable function on the bounded interval.} We then sample these curves using a uniform partition on $t\in[-0.5,0.5]$ with grid size $\Delta t$ ranging from $1/16$ to $\correct{1/4096}$ and project these data points in $\mathbb{R}^3$ onto $\mathbb{S}^2$ using the normalization $\mathbf{y}_i=\mathbf{x}_i/\|\mathbf{x}_i\|$ with \correct{$\mathbf{x}_i=(1,y(t_i),z(t_i))$}. We denote the interpolant as $\mathbf{y}(t)$ and the exact projection of the underlying curve onto the sphere as $\mathbf{z}(t)$ for $t\in[-0.5,0.5]$. Then we define the error between the reconstruction and the curve generated by the underlying function {$\mathbf{z}(t)$} by 
$$
e_{\mathbf{y}}=\int_{-0.5}^{0.5} \| \mathbf{y}(t)-\mathbf{z}(t) \| \, dt \, .
$$
Numerically, we approximate the integral using the triangles with vertices taken from $\mathbf{y}(t)$ and $\mathbf{z}(t)$. Note that since both $\mathbf{y}(t)$ and $\mathbf{z}(t)$ give points on a sphere, the definition of such an error does not provide the area of the mismatch on $\mathbb{S}^2$ but only certain cross-sessions in the Cartesian space. Finally, once we have obtained the error from each set of sampling points, we determine the numerical rate of convergence using $\rho_{\Delta t}=\log_2(e_{\Delta t}/e_{\Delta t/2})$. 
\correct{For the non-differentiable case, we repeat the above process but replace the generating function by $z(t)=|f(t)|$. The absolute value will create a kink at $t=0$ while maintaining the smoothness of the curve elsewhere on $t\in(-0.5,0.5)$. Table \ref{tab:area_err} shows the error and the estimated order of SIDER, SQUAD, SENO2 and SENO3 when the underly generating function is non-differentiable and smooth, respectively.}

Here are our observations about the numerical convergence for all examples using SLERP, SQUAD, SENO2, and SENO3. 

\begin{itemize}
\item The piecewise linear interpolation SLERP achieves second-order convergence. Both errors and convergence orders are independent of the smoothness of the underlying curve. This observation is expected since the reconstruction uses only adjacent data points and is not affected by any kink in the given data.
\item The interpolation SQUAD uses 4 data points for local reconstructions. If the underlying curve is smooth enough, we see a third-order convergence. When the underlying curve is only $C^0$ where there is a kink, the convergence rate drops from three to two. \correct{We have bolded these numbers in Table \ref{tab:area_err}.}
\item Because of the ENO idea that we locally determine a set of stencils producing the least variations, our proposed SENO can avoid interpolation across the kink. The numerical accuracy of SENO2 and SENO3 are consistently three and four, respectively, independent of the smoothness of the underlying curve. 
\item In terms of the data required, each SQUAD interpolation takes 4 data points, while a SIDER2 (and therefore the final SENO2) interpolation needs only 3 data points. Still, both interpolation methods achieve the same order of convergence for a smooth enough generating function. However, when there is a kink in the underlying curve, SENO2 performs significantly better than SQUAD in both the numerical accuracy and the convergence rate, as shown in Table \ref{tab:area_err}(a).
\end{itemize}

\begin{figure}[!ht]
\centering
(a)\includegraphics[width=0.46\textwidth,clip]{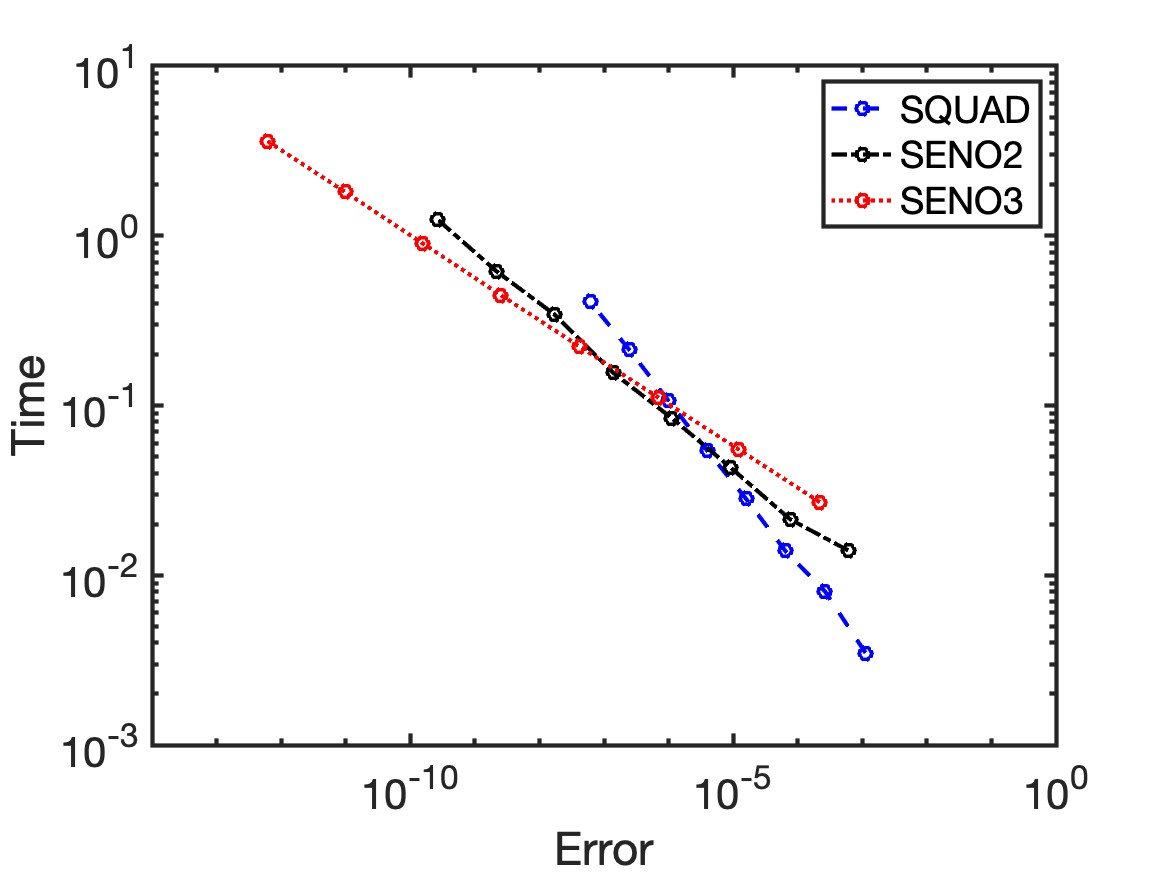}
(b)\includegraphics[width=0.46\textwidth,clip]{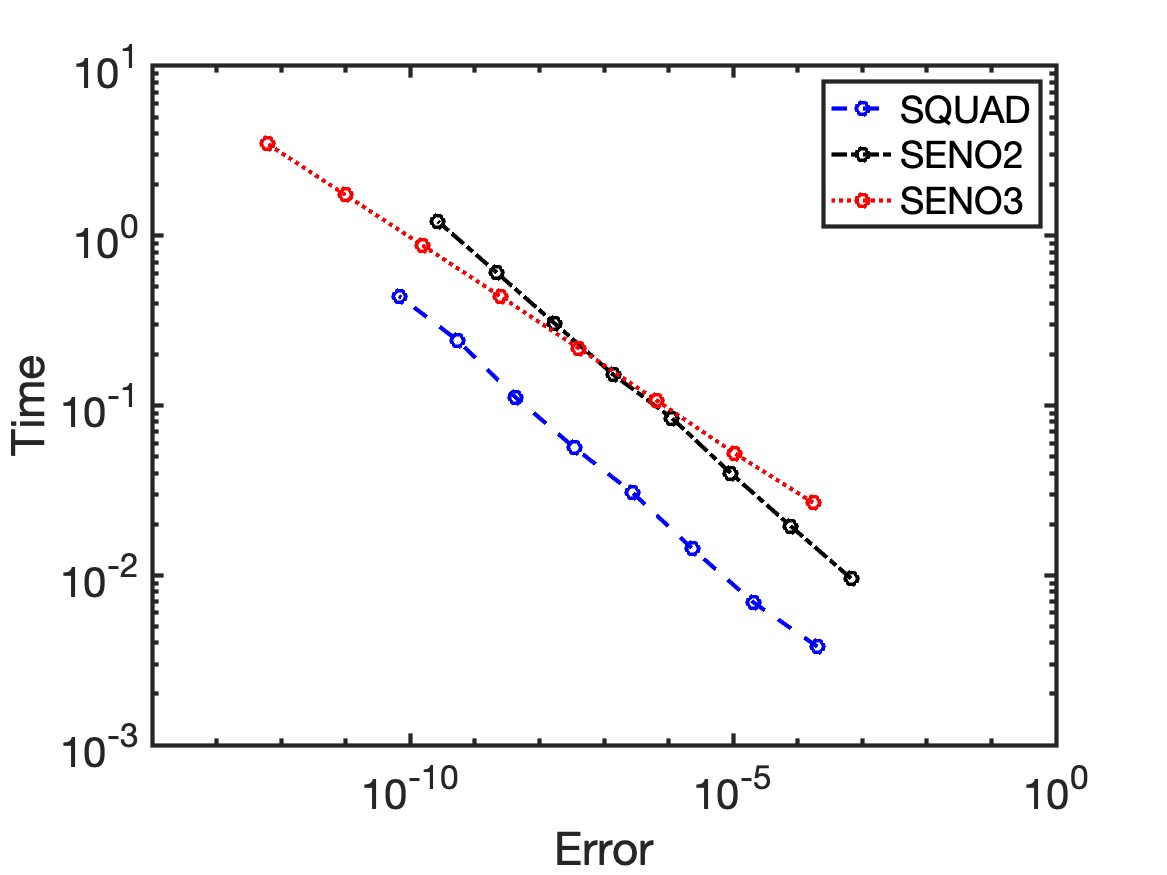}
\caption{\correct{(Section \ref{Sec:Efficiency}) The computational time versus the numerical error in the solution obtained by SLERP, SQUAD, SENO2 and SENO3. The underlying curve is (a) non-differentiable and (b) smooth.}}
\label{Fig:TimeVSError}
\end{figure}

\subsection{\correct{Efficiency}}
\label{Sec:Efficiency}

\correct{We have provided some complexity analysis in Table \ref{Table:Comparison} and have concluded that SENO might seem computationally less efficient than SQUAD since SENO3 (for example) takes 21 calls of SLERP while SQUAD requires only 3. Following the previous example, we are going to have more careful studies on the efficiency. We want to investigate if the increase in the reconstruction order is worth the computational complexity. Following the same study as in Section \ref{Sec:AccuracyConvergence}, we also measure the CPU-time corresponding to each configuration. Figure \ref{Fig:TimeVSError} shows the plots of the computational time versus the error in the reconstruction using three different interpolation methods. The curves for both SENO2 and SENO3 are almost identical for the non-differentiable and smooth cases.}

\correct{When the generating function is smooth, we see from Figure \ref{Fig:TimeVSError}(b) that SQUAD seems more efficient for a wide range of sampling densities. To achieve an accuracy of at least $O(10^{-10})$, for example, SQUAD interpolation takes less computational time than SENO3. Even though the convergence order of SQUAD is smaller than SENO3, it is affordable to refine the mesh to reduce the error in the interpolation. If the required error in the interpolation is further reduced, we might see that the blue dashed line and the red dotted line intersects, indicating that it is more efficient to use SENO3 than SQUAD. However, the accuracy level might be too close to machine epsilon to observe the benefit.}

\correct{However, when the underlying function is non-differentiable, we see from Figure \ref{Fig:TimeVSError}(a) that the curve corresponding to SENO3 is on the bottom compared to SQUAD and SENO2. This property indicates that the SENO3 is the most efficient approach when we want an accurate interpolant. To achieve a relatively small error, i.e. drawing a vertical line on the left regime, the computational time for SENO3 (the red dotted line) is the smallest. This figure clearly demonstrates the importance of the SENO reconstruction proposed in this work. When the data contains sharp turns due to a possible kink in the underlying curve, SENO can avoid interpolation across the singularity in the geometry and produce a high-order reconstruction.}

\begin{figure}[!htb]
\centering
(a)\includegraphics[width=0.46\textwidth]{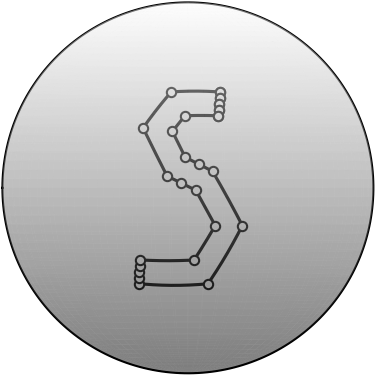}
(b)\includegraphics[width=0.46\textwidth]{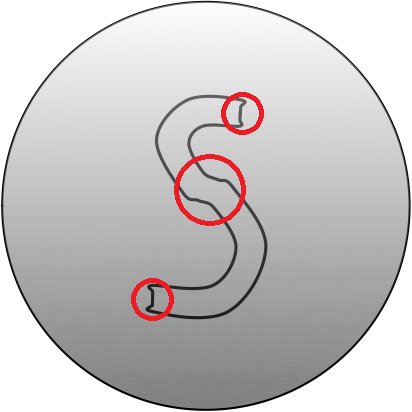}
(c)\includegraphics[width=0.46\textwidth]{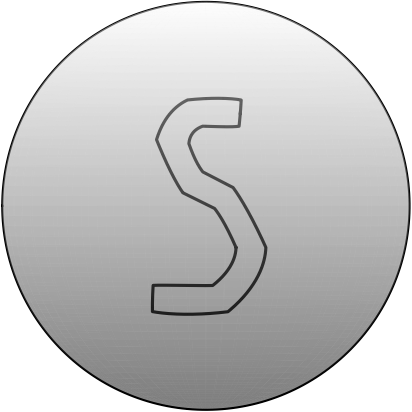}
(d)\includegraphics[width=0.46\textwidth]{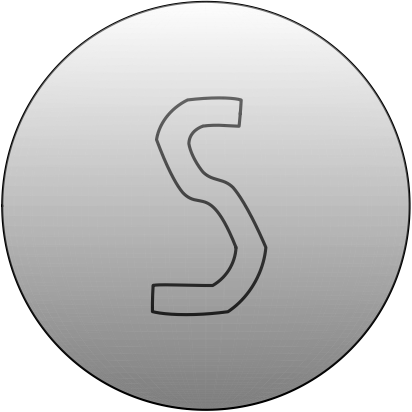}
\caption{(Section \ref{SubSec:LetterS}) The interpolation results about the data points on the unit sphere that resemble the letter \textit{S} using (a) SLERP, (b) SQUAD, (c) SENO2 and (d) SENO3. We have highlighted the oscillations obtained by SQUAD in (b).} 
\label{fig:ust_select}
\end{figure}

\begin{figure}[!htb]
\centering
(a)\includegraphics[width=0.31\textwidth]{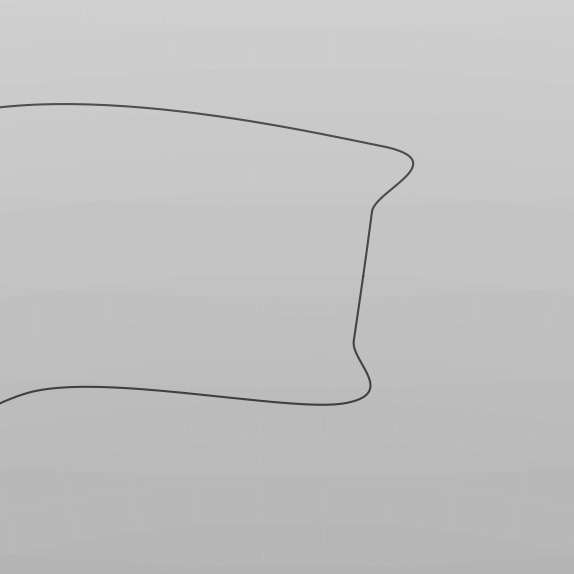}
\includegraphics[width=0.31\textwidth]{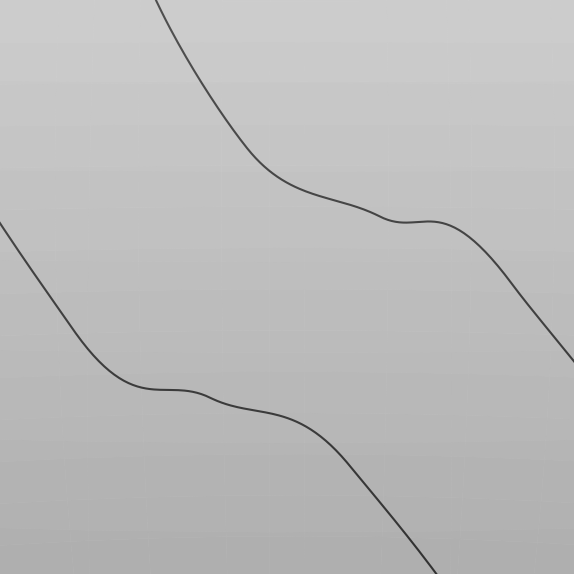}
\includegraphics[width=0.31\textwidth]{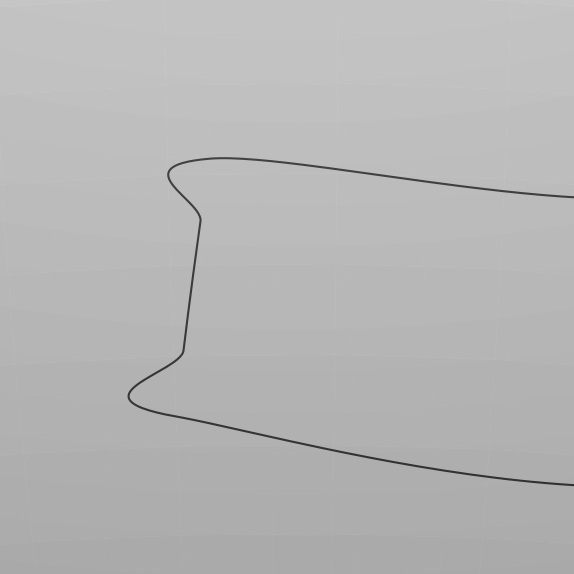} \\
(b)\includegraphics[width=0.31\textwidth]{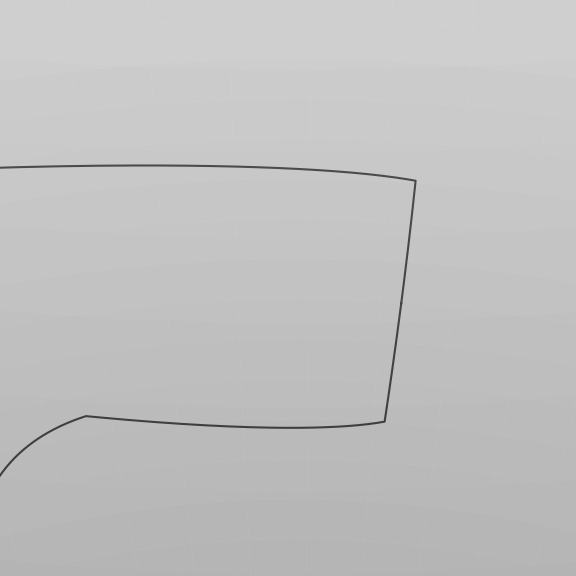}
\includegraphics[width=0.31\textwidth]{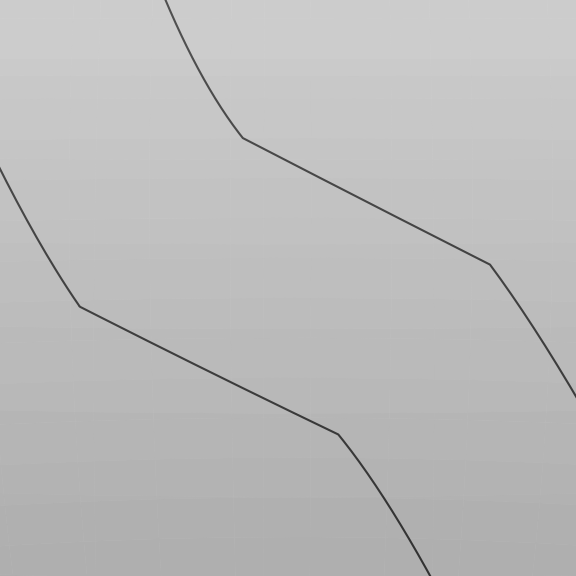}
\includegraphics[width=0.31\textwidth]{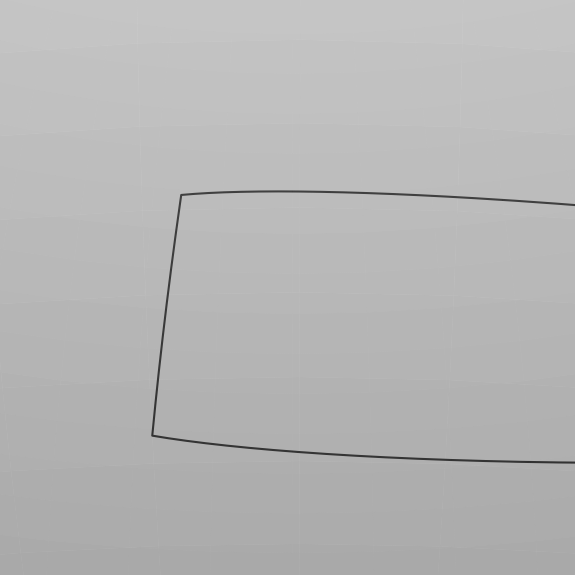} \\
(c)\includegraphics[width=0.31\textwidth]{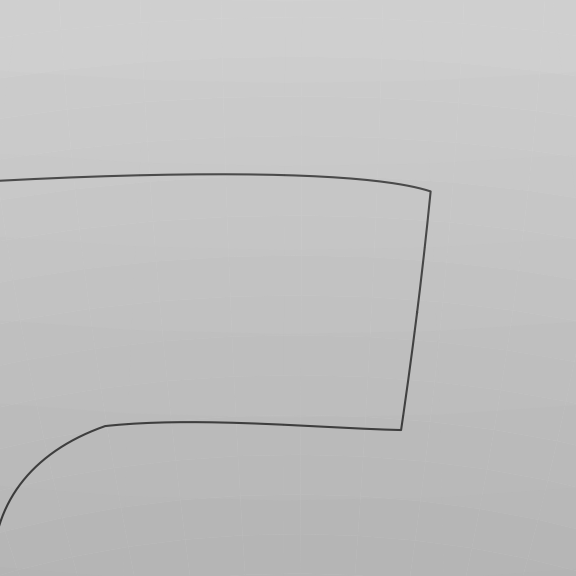}
\includegraphics[width=0.31\textwidth]{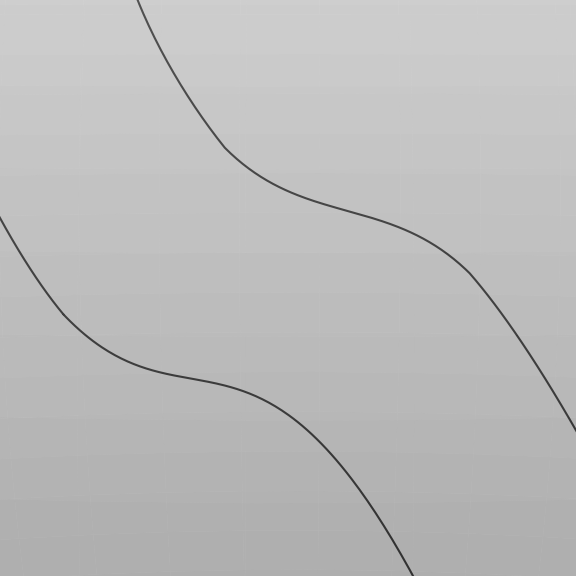}
\includegraphics[width=0.31\textwidth]{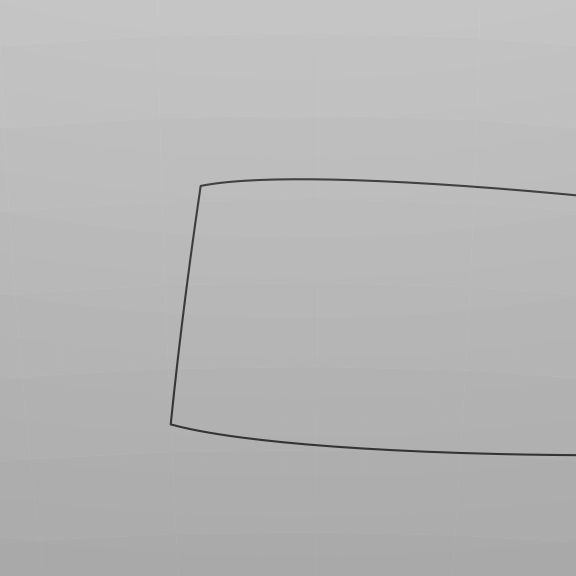}
\caption{(Section \ref{SubSec:LetterS}) Zoom-in of various regions of (a) the SQUAD, (b) the SENO2 and (c) the SENO3 interpolants in Figure \ref{fig:ust_select}.}
\label{fig:s_squad_closeup}
\end{figure}

\subsection{The Letter \textit{S} on Sphere}
\label{SubSec:LetterS}

Figure \ref{fig:ust_select} demonstrates how various interpolation schemes mentioned in this paper interpolate the data points on $\mathbb{S}^2$ with data points sampled from the boundary of the letter \textit{S}. It is clear that SLERP only creates geodesics between the neighboring data points, while SQUAD produces oscillations near the middle part and the cusps around the corners of the character, as highlighted in Figure \ref{fig:ust_select}(b). A zoom-in of these parts are shown in Figure \ref{fig:s_squad_closeup}. On the other hand, both SIDER2 and SIDER3 can control these overshoots well. We see sharp corners near these kinks and smooth interpolants when the segment is smooth, and we have an adequate number of sample points.

\section{Conclusion}

This article develops high-order interpolation schemes for data points on $\mathbb{S}^2$. We present the SIDER interpolation constructed based on an approach similar to the B\'{e}zier curves. To improve the accuracy of the reconstruction when the underlying function is not smooth enough, we also propose a SENO reconstruction by extending the ENO idea for data points in $\mathbb{R}^n$ to {$\mathbb{S}^2$}. \correct{A possible future work is to follow a similar strategy as in the weighted-ENO (WENO) approach \cite{loc94,shu97,jiapen00} to combine multiple SIDER's and achieve a higher-order reconstruction with weighting determined by the smoothness of the underlying curve.}

\section*{Acknowledgment}
The work of Leung was supported in part by the Hong Kong RGC grant 16302819. 

\bibliographystyle{plain}
\bibliography{syleung,referencefile}


\section*{Appendix: Derivatives of SQUAD, SIDER2 and SIDER3}

\subsection*{Time Derivatives of the Interpolants}

Assume $t \in [0, 1]$, and $\mathbf{q_i}=(0,\mathbf{p_i})$ is the starting point. The derivatives of  SLERP are given by
\begin{equation*}
    \dfrac{d}{dt}(\text{SLERP}(\mathbf{q_i}, \mathbf{q_{i+1}}, f(t)))^{\pm 1} = \pm \text{SLERP}(\mathbf{q_i}, \mathbf{q_{i+1}}, f(t)) \cdot \ln((\mathbf{q_i})^{-1} \mathbf{q_{i+1}}) \cdot \dfrac{df(t)}{dt}.
\end{equation*} 

To simplify the expressions, we introduce the following \Yuki{notations}.
\begin{equation*}
    \begin{array}{cccc}
        \text{Scheme} & \text{SIDER2} & \text{SQUAD} & \text{SIDER3} \\
        f(t) & t & 2t(1-t) & t \\
        g(t) & t & t & 3t/2 \\
        h(t) & t & t & (3t-1)/2 \\
        \mathbf{p}(t) & \text{SLERP}(\mathbf{q_i}, \mathbf{d_{(i+1)a}}, t) & \text{SLERP}(\mathbf{q_i}, \mathbf{q_{i+1}}, t) & \text{SIDER2}(\mathbf{q_i}, \mathbf{q_{i+1}}, \mathbf{q_{i+2}}, t) \\
        \mathbf{s}(t) & \text{SLERP}(\mathbf{d_{(i+1)b}}, \mathbf{q_{i+2}}, t) & \text{SLERP}(\mathbf{s_i}, \mathbf{s_{i+1}}, t) & \text{SIDER2}(\mathbf{q_{i+1}}, \mathbf{q_{i+2}}, \mathbf{q_{i+3}}, t) \\
\end{array}
\end{equation*} 

Since 
\begin{equation*}
    \begin{split}
        &\text{{ }{ }{ }{ }} (0, \mathbf{p}(g(t))) \times ((0, -\mathbf{p}(g(t))) \times (0, \mathbf{s}(h(t))))^{f(t)} \\
        &= (0, \mathbf{p_g}(t)) \times ((0, -\mathbf{p_g}(t)) \times (0, \mathbf{s_h}(t))^{f(t)}) = (0, \mathbf{p_g}(t)) \times (\mathbf{p_g}(t) \cdot \mathbf{s_h}(t), -\mathbf{p_g}(t) \times \mathbf{s_h}(t))^{f(t)} \\
        &= (0, \mathbf{p_g}(t)) \times (\cos(\theta_{\mathbf{ps}}(t)), \sin(\theta_{\mathbf{ps}}(t)) \mathbf{a_{ps}}(t))^{f(t)} = (0, \mathbf{p_g}(t)) \times (\cos(f_{\theta, \mathbf{ps}}(t)), \sin(f_{\theta, \mathbf{ps}}(t)) \mathbf{a_{ps}}(t)),
    \end{split}
\end{equation*} 
the first and the second time derivatives are given by
\begin{eqnarray*} 
&&\dfrac{d}{dt} \left[(0, \mathbf{p}(g(t))) \times ((0, \mathbf{p}(g(t))) \times (0, \mathbf{s}(h(t))))^{f(t)}\right] \\
&=& (0, (\mathbf{p_g})'(t)) \times (\cos(f_{\theta, \mathbf{ps}}(t)), \sin(f_{\theta, \mathbf{ps}}(t)) \mathbf{a_{ps}}(t)) - \mathbf{p_g}(t) \cdot \sin(f_{\theta, \mathbf{ps}}(t)) \cdot (f_{\theta, \mathbf{ps}})'(t) \\
& & + (0, \mathbf{p_g}(t)) \times (0, \cos(f_{\theta, \mathbf{ps}}(t)) \cdot (f_{\theta, \mathbf{ps}})'(t) \times \mathbf{a_{ps}}(t) + \sin(f_{\theta, \mathbf{ps}}(t)) \times (\mathbf{a_{ps}})'(t))
\end{eqnarray*} 
and
\begin{eqnarray*} 
&& \dfrac{d}{dt} \left[\dfrac{d}{dt} \left[(0, \mathbf{p}(g(t))) \times ((0, -\mathbf{p}(g(t))) \times (0, \mathbf{s}(h(t))))^{f(t)}\right]\right] \\
&=& (0, -(\theta_{\mathbf{p}})^2(\mathbf{p_g}(t))) \times (\cos(f_{\theta, \mathbf{ps}}(t)), \sin(f_{\theta, \mathbf{ps}}(t)) \mathbf{a_{ps}}(t)) - 2(\mathbf{p_g})'(t) \cdot \sin(f_{\theta, \mathbf{ps}}(t)) \cdot (f_{\theta, \mathbf{ps}})'(t) \\
&&+ 2 \cdot (0, (\mathbf{p_g})'(t)) \times (0, \cos(f_{\theta, \mathbf{ps}}(t)) \cdot (f_{\theta, \mathbf{ps}})'(t) \times \mathbf{a_{ps}}(t) + \sin(f_{\theta, \mathbf{ps}}(t)) \times (\mathbf{a_{ps}})'(t)) \\
&&- \mathbf{p_g}(t) \cdot (\cos(f_{\theta, \mathbf{ps}}(t)) \cdot ((f_{\theta, \mathbf{ps}})'(t))^2 + \sin(f_{\theta, \mathbf{ps}}(t)) \cdot (f_{\theta, \mathbf{ps}})''(t)) \\
&&+ (0, \mathbf{p_g}(t)) \times (0, (\cos(f_{\theta, \mathbf{ps}}(t)) \cdot (f_{\theta, \mathbf{ps}})''(t) - \sin(f_{\theta, \mathbf{ps}}(t)) \cdot ((f_{\theta, \mathbf{ps}})'(t))^2) \times \mathbf{a_{ps}}(t)) \\
&&+ (0, \mathbf{p_g}(t)) \times (0, 2 \cos(f_{\theta, \mathbf{ps}}(t)) \cdot (f_{\theta, \mathbf{ps}})'(t) \times (\mathbf{a_{ps}})'(t) + \sin(f_{\theta, \mathbf{ps}}(t)) \times (\mathbf{a_{ps}})''(t))  
\end{eqnarray*} 
where 
\begin{eqnarray*} 
\theta_{\mathbf{p}} &=& \arccos(\mathbf{p_i} \cdot \mathbf{c_{(i+1)a}} \text{ for SIDER2 or } \cdot \mathbf{p_{i+1}} \text{ for SQUAD}), \\
\mathbf{a_p} &=& (-\mathbf{p_i} \times \mathbf{c_{(i+1)a}} \text{ for SIDER2 or } \times \mathbf{p_{i+1}} \text{ for SQUAD}) / \sin(\theta_{\mathbf{p}}), \\
(0, \mathbf{p_g}(t)) &=& (0, \mathbf{p}(g(t))) = (\cos(g(t) \cdot \theta_{\mathbf{p}}))\mathbf{p_i} + (\sin(g(t) \cdot \theta_{\mathbf{p}}))(\mathbf{a_p} \times \mathbf{p_i}), \notag \\
(0, \mathbf{s_h}(t)) &=& (0, \mathbf{s}(h(t))), \\
(\mathbf{p_g})'(t) &=& (\theta_{\mathbf{p}})(\mathbf{a_p} \times \mathbf{p_g}(t)), (\mathbf{p_g})''(t) = -(\theta_{\mathbf{p}})^2(\mathbf{p_g}(t)) \\
\theta_{\mathbf{ps}}(t) &=& \arccos(\mathbf{p_g}(t) \cdot \mathbf{s_h}(t)), \\
(\theta_{\mathbf{ps}})'(t) &=& -((\mathbf{p_g})'(t) \cdot \mathbf{s_h}(t) + \mathbf{p_g}(t) \cdot (\mathbf{s_h})'(t)) / \sin(\theta_{\mathbf{ps}}(t)), \\
\Yuki{(\theta_{\mathbf{ps}})''(t)} &=& -((\mathbf{p_g})''(t) \cdot \mathbf{s_h}(t) + 2 \cdot (\mathbf{p_g})'(t) \cdot (\mathbf{s_h})'(t) + \mathbf{p_g}(t) \cdot (\mathbf{s_h})''(t)\\
&&+ \cos(\theta_{\mathbf{ps}}(t)) \cdot ((\theta_{\mathbf{ps}})'(t))^2) / \sin(\theta_{\mathbf{ps}}(t)), \notag \\
\mathbf{a_{ps}}(t) &= &(-\mathbf{p_g}(t) \times \mathbf{s_h}(t))/\sin(\theta_{\mathbf{ps}}(t)), \\
(\mathbf{a_{ps}})'(t) &=& (-(\mathbf{p_g})'(t) \times \mathbf{s_h}(t) + (-\mathbf{p_g}(t)) \times (\mathbf{s_h})'(t) \\ 
&&- \cos(\theta_{\mathbf{ps}}(t)) \times (\theta_{\mathbf{ps}})'(t) \times \mathbf{a_{ps}}(t))/\sin(\theta_{\mathbf{ps}}(t)), \notag \\
(\mathbf{a_{ps}})''(t) &=& (-(\mathbf{p_g})''(t) \times \mathbf{s_h}(t) + 2 \cdot (-(\mathbf{p_g})'(t)) \times (\mathbf{s_h})'(t) + (-\mathbf{p_g}(t)) \times (\mathbf{s_h})''(t) \notag \\
&&+ \sin(\theta_{\mathbf{ps}}(t)) \times ((\theta_{\mathbf{ps}})'(t))^2 \times \mathbf{a_{ps}}(t) - \cos(\theta_{\mathbf{ps}}(t)) \times (\theta_{\mathbf{ps}})''(t) \times \mathbf{a_{ps}}(t) \\
&&- 2 \cdot \cos(\theta_{\mathbf{ps}}(t)) \times (\theta_{\mathbf{ps}})'(t) \times (\mathbf{a_{ps}})'(t))/\sin(\theta_{\mathbf{ps}}(t)), \notag \\
f_{\theta, \mathbf{ps}}(t) &=& f(t) \cdot \theta_{\mathbf{ps}}(t), \quad \dfrac{d^n f_{\theta, \mathbf{ps}}(t)}{dt^n} = \sum_{r = 0}^n \dfrac{d^r f(t)}{dt^r} \cdot \dfrac{d^{n-r} \theta_{\mathbf{ps}}(t)}{dt^{n-r}}. \text{ (product rule}) \, .
\end{eqnarray*} 
The derivatives of $\mathbf{s_h}(t)$ follows the same manner as $\mathbf{p_g}(t)$, but they are not explicitly written out because $\mathbf{p}(t), \mathbf{s}(t), f(t), g(t)$ and/or $h(t)$ of SQUAD, SIDER2 and SIDER3 are not the same.

\subsection*{Angular Derivatives of the Interpolants}

With references to \cite{dantam14}, if the interpolant is $\mathbf{q}(t)$, then
we can deduce that
\begin{equation*}
    \begin{split}
        \mathbf{\omega}(t) &= 2 \mathbf{q}'(t) \times (\mathbf{q}(t))^{-1}, \\
        \mathbf{\alpha}(t) &= (2\mathbf{q}''(t) - \mathbf{\omega}(t) \times \mathbf{q}'(t)) \times (\mathbf{q}(t))^{-1}, \text{ and} \\
        \mathbf{\zeta}(t) &= (2\mathbf{q}'''(t) - 2 \mathbf{\alpha}(t) \times \mathbf{q}'(t) - \mathbf{\omega}(t) \times \mathbf{q}''(t)) \times (\mathbf{q}(t))^{-1}.
    \end{split}
\end{equation*}

\end{document}